\documentclass[12pt]{amsart}
\usepackage{amsfonts,amssymb}
\input xy
\xyoption{all}

\newtheorem{theo}{Theorem}[section]
\newtheorem{pro}[theo]{Proposition}
\newtheorem{lem}[theo]{Lemma}

\newtheorem{cor}[theo]{Corollary}
\newtheorem{defin}[theo]{Definition}
\newcommand{\ra}{\rightarrow}

\newcommand{\Gr}{{\rm Gr}}
\newcommand{\N}{{\mathbb N}}
\newcommand{\Q}{{\mathbb Q}}
\newcommand{\R}{{\mathbb R}}
\newcommand{\bfsigma}{\sigma\hspace{-0.24cm}\sigma\hspace{-0.24cm}\sigma}
\newcommand{\LolocG}{L^1_{\rm loc}(G)}

\title{Connected Lie groups and property RD}
\author{I. Chatterji$^\dagger$, Ch. Pittet and L. Saloff-Coste$^{\dagger\dagger}$}

\thanks{$\dagger$ Partially supported by the Swiss Science Foundation grant PA002-101406 and by NSF grant DMS 0405032.}
\thanks{${\dagger\dagger}$ Partially supported by NSF grant DMS 0102126.}
\thanks{2000 Mathematics Subject Classification number: 22D15, 22E30, 43A15 and  46L05.}
\date{\today}
\begin{document}
\maketitle
\begin{abstract} For a locally compact group, property RD gives a control on the convolution norm of any compactly supported measure in terms of the $L^2$-norm of its density and the diameter of its support. We give a complete classification of those Lie groups with property RD.
\end{abstract}
%%%%%%%%%%%%%%%%%%%%%%%%%
\section*{Introduction}
%%%%%%%%%%%%%%%%%%%%%%%
The property of Rapid Decay (property RD) emerged from the work of Haagerup in \cite{Haa} and was first studied systematically by P. Jolissaint in \cite{Jol} who mostly worked in the context of finitely generated groups. For a locally compact group, property RD gives a control on the convolution norm (i.e., operator norm) of any compactly supported measure in terms of the $L^2$-norm of its density and the diameter of its support. The terminology ``Rapid Decay''comes from the fact that a group has property RD if and only if any ``rapidly decaying'' function  is an $L^2$ convolutor (see Definition \ref{rdfunct} and Lemma \ref{EquivRD} below). Property RD is relevant in the context of C*-algebras. In particular V. Lafforgue used it in \cite{Laff} to prove the Baum-Connes conjecture for some groups having property~(T). The authors observed \cite{heatker} that property RD is also relevant to the study of the asymptotic behavior of the probability of return of random walks on non-amenable groups (for a general reference on random walks, see \cite{Woess}). This will be used here in Section \ref{class} to make the link between property RD and Varopoulos' work \cite{V}. 

The main result of this paper is a precise algebraic description of those connected Lie groups that have property RD (in what follows, Lie groups are real Lie groups).
\begin{theo}\label{classe} Let $G$ be a connected Lie group. The following are equivalent.
\begin{itemize}
\item[(a)] $G$ has property RD.
\item[(b)] The Lie algebra ${\mathbf g}$ of $G$ decomposes as a direct product ${\mathbf g}={\mathbf s}\times{\mathbf q}$, where ${\mathbf s}$ is semisimple or $\{0\}$ and ${\mathbf q}$ is an algebra of type R.
\item[(c)] The universal cover $\widetilde{G}$ of $G$ decomposes as a direct product $\widetilde{S}\times\widetilde{Q}$, where $\widetilde{S}$ is semisimple and $\widetilde{Q}$ has polynomial volume growth.\end{itemize}\end{theo}
We will also extend this result to compactly generated virtually connected groups in Corollary \ref{MerciBekka}. The equivalence between (b) and (c) is well-known (see e.g. \cite{Var}, \cite{G}, \cite{Jen}). That (a) implies (b) follows from Varopoulos' work in \cite{V}; this will be explained in Section \ref{class}. That (c) implies (a) occupies a large portion of this work. A short description of the paper is as follows. Notations are set in Section \ref{basics}. Sections \ref{RD}, \ref{stabilitefacile} and \ref{amen} discuss property RD in the context of locally compact groups (see also \cite{Jol} and \cite{JiSch}). In Section \ref{NAaRDrad} we show that certain solvable groups that do not have property RD nevertheless satisfy a weaker radial version. Namely, on these groups, compactly supported measures with radial densities satisfy the RD inequality. In Section \ref{ss} we use the results of Section \ref{NAaRDrad} to establish property RD for semisimple Lie groups with finite center. This was proved independently by V. Lafforgue \cite{LV} in an unpublished manuscript and was understood by others including M. Cowling or N. Higson, using \cite{CHH}. Our proof differs from Lafforgue's \cite{LV} and makes crucial use of a lemma from \cite{CGHM}. In Section \ref{ext} we establish the stability of property RD under some central extensions, a result proved by Jolissaint \cite{Jol} in the case of finitely generated discrete groups. In Section \ref{gss} we show that Theorem \ref{classe}(c) implies Theorem \ref{classe}(a). Section \ref{class} concludes the proof of Theorem \ref{classe}. 
%%%%%%%%%%%%%%%%%%%%%%%%%%%%
\section{Basic notation}\label{basics}
%%%%%%%%%%%%%%%%%%%%%%%%%%%%%%
Throughout this paper all measures are Borel regular and all groups are locally compact groups. Let $G$ be such a group. For $f$ a continuous function on $G$, set $\check{f}(x)=\overline{f(x^{-1})}$. For a Borel measure $\mu$, define the measure $\check{\mu}$ by $\check{\mu}(f)=\overline{\mu(\check{f})}$. Given two finite complex Borel measures $\mu_1,\mu_2$, the convolution $\mu_1*\mu_2$ is a finite complex Borel measure defined by
$$\mu_1*\mu_2(f)=\int_{G\times G} f(xy)d\mu_1(x)d\mu_2(y),$$
for any $f\in{\mathcal C}_0(G)$, where ${\mathcal C}_0(G)$ denotes continuous compactly supported functions on $G$. Let us denote by $\vec{\nu}(dx)=\vec{d}x$ a left Haar measure, so that $\nu(dx)=dx=\vec{\nu}\check{\phantom v}(dx)$ is a right Haar measure (left and right Haar measures are unique up to a multiplicative constant). Denote by $L^2(G)$ the Hilbert space $L^2(G,\vec{\nu})$ equipped with the inner product
$$\left<f,g\right>=\int_Gf(x)\overline{g(x)}\vec{d}x.$$
The modular function $m$, defined by $\vec{\nu}(Vg)=m(g)\vec{\nu}(V)$ (for any Borel set $V$ and $g\in G$) relates left and right measures as follows
$$dx= m(x^{-1})\vec{d}x=m(x)^{-1}\vec{d}x.$$
For a function $f\in {\mathcal C}_0(G)$ and a finite complex Borel measure $\mu$, we set
\begin{eqnarray*}L_\mu f(x)& = & \mu*f(x)=\int_Gf(y^{-1}x)d\mu(y)\\
R_\mu f(x)& =&f*\mu(x)=\int_Gm(y)^{-1}f(xy^{-1})d\mu(y).\end{eqnarray*}
Given two functions $f_1,f_2\in {\mathcal C}_0(G)$, the convolution $f_1*f_2$ 
is the function
$$f_1*f_2(x)=\int_G f_1(xy)f_2(y^{-1})\vec{d}y=\int_G f_1(y) f_2(y^{-1}x)\vec{d}y.$$
These definitions are consistent when the measure $\mu_1$ (respectively $\mu_2$) has density $f_1$ (respectively $f_2$) with respect to the left Haar measure $\vec{\nu}$ (see Chapter V of \cite{HR}) in the sense that
\begin{eqnarray*}\mu_1*f&=&f_1*f,\\
 f*\mu_2&=&f*f_2\end{eqnarray*}
for any $f\in{\mathcal C}_0(G)$ and the density of $\mu_1*\mu_2$ is $f_1*f_2$. We refer to \cite[Paragraph 20]{HR} for background information. For a Borel measure $\mu$ with finite total mass $\|\mu\|=\int_Gd|\mu|(y)$, $L_\mu$ is a continuous operator on $L^2(G)$ and $\|L_\mu\|_{2\ra 2}\le \|\mu\|$, where $\|L_\mu\|_{2\ra 2}$ denotes the operator norm of $L_\mu$ as an operator acting on $L^2(G)$, namely
$$\|L_\mu\|_{2\ra 2}=\sup_{f\in L^2(G)}\frac{\|L_{\mu}(f)\|_2}{\|f\|_2}.$$
For a Borel measure $\mu$ such that $\int_G m(y)^{-1/2}d|\mu|(y)<\infty$, set $\hat{\mu}(dx)=m(x)^{-1/2}\mu(dx)$. Notice that the map
\begin{eqnarray*}I:L^2(G)&\to & L^2(G)\\
f &\mapsto & m(x)^{-1/2}\check{f}\end{eqnarray*}
is unitary and that $R_{\overline{\mu}}=I\circ L^*_{\hat{\mu}}\circ I$. Hence $R_{\mu}$ is a continuous operator on $L^2(G)$ and
\begin{equation}\label{def}\|R_\mu\|_{2\ra 2}=\|L_{\hat{\mu}}\|_{2\ra 2}\le \int_G m(y)^{-1/2}d|\mu|(y).\end{equation}
%%%%%%%%%%%%%%%%%%%%%%%
\section{Property RD}\label{RD}
%%%%%%%%%%%%%%%%%%%%%
There are two possible definitions of property RD on a general locally compact group. They differ a priori only if the group is not unimodular. We call them RD and RD'. We will show that in fact RD and RD' are equivalent. It has been shown by Ji and Schweitzer, Theorem 2.2 in \cite{JiSch}, that a topological locally compact group with property RD has to be unimodular, but our proof that RD is equivalent to RD' doesn't use this fact. We shall keep the distinction between RD and RD' because these properties will eventually differ on some subsets of $L^2(G)$ (such as radial functions, see Section \ref{NAaRDrad}). Recall that a length function on a locally compact group $G$ is a Borel map $L:G\to{\mathbb R}^+$ satisfying $L(1)=0$, $L(gh)\leq L(g)+L(h)$ and $L(g)=L(g^{-1})$. For any locally integrable non-negative function $f$ in $G$, we define the Borel measure $F$ by $F(V)=\int_Vf(x)\vec{d}x$ (where $V\subseteq G$ is compact). We say that $f\in\LolocG$ has compact support if $F$ (as a measure) has compact support. We shall denote by $B_L(R)$ the $L$-ball of radius $R\geq 1$ centered at the identity $1\in G$.
\begin{defin}\label{DefRDRD'}
Let $E\subseteq L^2(G)$, we say that the pair $(G,L)$ has \emph{property RD'$_E$} if there exists two constants $C,D\ge 0$ such that for any function $f\in E$ with compact support in $B_L(R)$, $R\geq 1$, we have
$$\|L_F\|_{2\ra 2}\le CR^D\|f\|_2.$$
We say that the pair $(G,L)$ has \emph{property RD$_E$} if there exists two constants $C,D\ge 0$ such that for any function $f\in E$ with compact support in $B_L(R)$, $R\geq 1$, we have
$$\|R_F\|_{2\ra 2}\le CR^D\|f\|_2.$$\end{defin}
In this paper, all topological groups will be compactly generated in the following sense. A topological group $G$ is compactly generated if it admits a compact neighborhood $K$ of the identity that generates $G$ (i.e. $K$ is such that $G=\bigcup_{n\in\N}K^n$). Note that a compactly generated group is then automatically locally compact. If $G$ is compactly generated and $K$ is a compact generating set satisfying $K=K^{-1}$, we call algebraic length function the length function defined by 
$$L_K(g)=\inf\{n: g\in K^n\},\ L_K(1)=0.$$ 
A length function $L$ on $G$ is said to be proper if it is bounded on any compact set $U$, namely 
$$M_U=\sup\{L(u)|u\in U\}<\infty.$$
It is straightforward to see that if $G$ is compactly generated then any algebraic length function is proper. Moreover, if $(G,L)$ has property RD for some proper length function $L$, then so will $(G,L_K)$ for any compact generating set $K$. Indeed, for $g=s_1\dots s_n$ with $s_i\in K$ and $n$ minimal:
$$L(g)=L(s_1\dots s_n)\le\sum_{i=1}^nL(s_i)\le M_KL_K(g).$$
Hence all algebraic length functions are comparable and we shall talk about algebraic length without specifying the compact generating set. Further examples of length functions can be obtained by letting the group $G$ act by isometries on a metric space $(X,d)$ and by setting $L(g)=d(x_0,g(x_0))$ for some base point $x_0\in X$.
\begin{defin}Let $E\subseteq L^2(G)$, we call $E$ \emph{sd-closed} if $E*\check{E}\subseteq E$, and \emph{sm-closed} if $E*m^{-1}\check{E}\subseteq E$.\end{defin}
The terminology sd and sm stands for symmetric density and symmetric measure. Indeed, $m^{-1}\check{E}$ is the set of densities of the measure $\check{F}$ where $F=f\vec{d}x$ with $f\in E$ and for instance the set $\{f\in L^2(G)|f=\check{f}\}$ is sd-closed whereas the set $\{f\in L^2(G)|\check{F}=F\}$ is sm-closed. The following lemma shows that for sd-closed (respectively sm-closed) sets it is enough to check the RD$_E$ inequality on measures with symmetric densities (respectively on symmetric measures).
\begin{lem}Let $G$ be a locally compact group and $L$ a proper length function on $G$.
\begin{itemize}
\item[(1)]Let $E\subseteq L^2(G)$ be sd-closed. If there are two constants $C,D\ge 0$ such that for any $f\in E$ with $f=\check{f}$ and supported on $B_L(R)$ for $R\geq 1$
$$\|R_F\|_{2\ra 2}\leq CR^D\|f\|_2,$$
then $(G,L)$ has property RD$_E$.
\item[(2)]Let $E\subseteq L^2(G)$ be sm-closed. If there are two constants $C,D\ge 0$ such that for any $f\in E$ with $f=m^{-1}\check{f}$ and supported on $B_L(R)$ for $R\geq 1$
$$\|L_F\|_{2\ra 2}\leq CR^D\|f\|_2,$$
then $(G,L)$ has property RD'$_E$.\end{itemize}\end{lem}
\begin{proof} Let us start by quoting from \cite{HR} Discussion (20.20). If a measure $F$ has density $f$, then the adjoint of the operator $L_F$ acting on $L^2(G)$ satisfies
\begin{equation}\label{adjointL}L_F^*=L_{\check{F}}\end{equation}
where the measure ${\check{F}}$ has density $m^{-1}\check{f}$. A similar computation shows that
\begin{equation}\label{Ladjoint}R_F^*=R_{m\check{F}}\end{equation}
so that the measure $m\check{F}$ has density $\check{f}$. Now, $R_F^*R_F=R_{F*m\check{F}}$ and $F*m\check{F}$ has symmetric density $f*\check{f}$. Assuming (1) we compute:
\begin{eqnarray*}\|R_F^*\|_{2\ra 2}^2&=&\|R_F^*R_F\|_{2\ra 2}\leq C(2R)^D\|f*\check{f}\|_2\\
&=&C(2R)^D\|R_F^*(f)\|_2\leq C(2R)^D\|R_F^*\|_{2\ra 2}\|f\|_2.\end{eqnarray*}
The first inequality comes from the facts that we assumed $E$ sd-closed, and that if $f$ is supported on $B_L(R)$, then $f*\check{f}$ is supported on $B_L(2R)$. Since $\|R^*_F\|_{2\ra 2}=\|R_F\|_{2\ra 2}$ we deduce property RD$_E$. Assuming (2) and since $L_F^*L_F=L_{\check{F}*F}$ the computation is similar:
\begin{eqnarray*}\|L_F^*\|_{2\ra 2}^2&=&\|L_F^*L_F\|_{2\ra 2}\leq C(2R)^D\|m^{-1}\check{f}*f\|_2\\
&=&C(2R)^D\|L_F^*(f)\|_2\leq C(2R)^D\|L_F^*\|_{2\ra 2}\|f\|_2.\end{eqnarray*}
Since $\|L^*_F\|_{2\ra 2}=\|L_F\|_{2\ra 2}$ we deduce property RD'$_E$.\end{proof}
\begin{cor}\label{RD<=>RD'}Let $L$ be a length function on a locally compact group $G$ and $E\subseteq L^2(G)$ which is both sd-closed and sm-closed, then $(G,L)$ has property RD$_E$ if and only if $(G,L)$ has property RD'$_E$.\end{cor}
\begin{proof}First notice that for $f=\check{f}$ a symmetric function, then 
$$\|f\|_2=\sqrt{\int_G|f(x^{-1})|^2m(x)^{-1}\vec{d}x}=\|m^{-1/2}f\|_2.$$
Now assume property RD'$_E$. For $f\in E$ symmetric supported on $B_L(R)$ we have
$$\|R_F\|_{2\ra 2}=\|L_{\hat{F}}\|_{2\ra 2}\leq CR^D\|m^{-1/2}f\|_2= CR^D\|f\|_2.$$
Conversely notice that if $f$ satisfies $f=m^{-1}\check{f}$, then
$$\|m^{-1/2}f\|_2=\sqrt{\int_G|f(x^{-1})|^2m(x)^{-2}\vec{d}x}=\|m^{-1}\check{f}\|_2=\|f\|_2,$$
so that assuming property RD$_E$, for $f=m^{-1}\check{f}$ supported on $B_L(R)$ we have
$$\|L_F\|_{2\ra 2}=\|R_{\hat{F}}\|_{2\ra 2}\leq CR^D\|m^{-1/2}f\|_2= CR^D\|f\|_2.$$\end{proof}
\begin{defin}Let $G$ be a compactly generated group, and $L$ be a proper length function. We say that $(G,L)$ \emph{has property RD} if it has RD$_E$ for $E=L^2(G)$. A compactly generated group $G$ \emph{has property RD} if the pair $(G,L)$ has property RD for one (hence any) algebraic length function $L$.\end{defin} 
Note that according to Corollary \ref{RD<=>RD'} we can check property RD by using either left or right convolution operators, which ever is more convenient. More generally, according to Ji and Schweitzer \cite{JiSch}, property RD implies unimodularity. In view of the discussion above, a compactly generated group has property RD as soon as there exists a proper length function $L$, algebraic or not, such that the pair $(G,L)$ has property RD. Note that on any non compact compactly generated group there exists proper length functions for which property RD fails. For instance, one can construct those using $L(g)=\log (1+L_K(g))$ where $L_K$ is an algebraic length function.
\begin{defin}\label{rdfunct} Let $L$ be a length function on a locally compact group $G$. For $k\geq 0$, define
$$H^k_L(G)=\{f\in L^2(G)|\int_G(1+L(x))^{2k}|f(x)|^2\vec{d}x<\infty\}$$
and $H^{\infty}_L(G)=\bigcap_{k\geq 0}H^k_L(G)$. The space $H^{\infty}_L(G)$ is called the \emph{space of rapidly decaying functions}, the decay at infinity is faster than any inverse of a polynomial in terms of the distance to the identity (or to any fixed base point).
\end{defin}
The space $H^{\infty}_L(G)\subseteq L^2(G)$ is a Fr\'echet space for the projective limit topology induced by the sequence of norms $\|f\|_{2,L,k}=\|(1+L)^kf\|_2$. Recall that the \emph{reduced $C^*$-algebra} of a locally compact group $G$ is the operator norm closure of compactly supported continuous functions over $G$, viewed as acting on $L^2(G)$ via the left regular representation (i.e. as $L_F$, where $dF=f\vec{d}x$, $f\in {\mathcal C}_0(G)$). Namely
$$C^*_r(G)=\overline{{\mathcal C}_0(G)}^{\|\ \|_{2\ra 2}}.$$
 The following lemma shows that Definition \ref{DefRDRD'} of property RD coincides with the one given by Jolissaint in \cite{Jol}, and used by Ji and Schweitzer in \cite{JiSch}.
\begin{lem}\label{EquivRD}Let $G$ be a locally compact group and $L$ a proper length function on $G$. Then the following are equivalent:
\begin{itemize}
\item[(1)]$(G,L)$ has property RD.
\item[(2)]$(G,L)$ has property RD$_E$ for $E={\mathcal C}_0(G)$.
\item[(3)]There is $k>0$ and $C>1$ such that, for any $f\in{\mathcal C}_0(G)$,
$$\|L_F\|_{2\ra 2}\leq C\|(1+L)^kf\|_2.$$
\item[(4)]$H^{\infty}_L(G)\subseteq C^*_r(G)$.
\end{itemize}\end{lem}
\begin{proof}That (1) implies (2) is obvious since ${\mathcal C}_0(G)\subseteq L^2(G)$. We look at (2) implies (3): Take $f\in{\mathcal C}_0(G)$ and write $f=\sum_{n=1}^{\infty}f_n$, where $f_n(x)=f(x)$ if $n-1\leq L(x) < n$, and $0$ otherwise. Then $\|f\|_2^2=\sum_{n=1}^{\infty}\|f_n\|_2^2$ and
\begin{eqnarray*}\|L_F\|_{2\ra 2} &\leq &\sum_{n=1}^{\infty}\|L_{F_n}\|_{2\ra 2}\leq C\sum_{n=1}^{\infty} n^D\|f_n\|_2=C\sum_{n=1}^{\infty}n^{-1} n^{D+1}\|f_n\|_2\\
&\leq &C\sqrt{\sum_{n=1}^{\infty}n^{-2}}\sqrt{\sum_{n=1}^{\infty}n^{2D+2}\|f_n\|_2^2}\leq C'\|(1+L)^{D+1}f\|_2,\end{eqnarray*}
for some finite constant $C'$. That (3) implies (4) is by density of ${\mathcal C}_0(G)$ in $H^{k}_L(G)$. Finally, let us show that (4) implies (1). First notice that the graph of the inclusion $H^{\infty}_L(G)\to C^*_r(G)$, $f\mapsto L_F$ is closed. Indeed, let $\{f_n\}_{n\in\N}$ in $H^{\infty}_L(G)$ tends to $f$ in $H^{\infty}_L(G)$ and $\{L_{F_n}\}_{n\in\N}$ tends to $g$ in $C^*_r(G)$. Let us prove that $L_F=g$. Obviously $\{\left<L_{F_n}\varphi,\psi\right>\}_{n\in\N}$ converges to $\left<g(\varphi),\psi\right>$ for any $\varphi,\psi\in{\mathcal C}_0(G)$. Since $\{f_n\}_{n\in\N}$ tends to $f$ in $L^2(G)$ as well, it implies that $\{f_n*\varphi*\psi(1)\}_{n\in\N}$ converges to $f*\varphi*\psi(1)$ for any $\varphi,\psi\in{\mathcal C}_0(G)$. As
$$f*\varphi*\psi(1)=\left<L_F\varphi,\check{\psi}\right>$$
we conclude that $L_F=g$. The closed graph theorem in the generality of Proposition 1, Chapter I page 20 of \cite{Bourb}, applied to the Fr\'echet space $H^{\infty}_L(G)$ and to $C^*_r(G)$ (viewed as a Banach space) then implies that the inclusion (4) is continuous. This by definition amounts to the existence of a $k>0$ and $C>1$ such that $\|L_F\|_{2\ra 2}\leq C\|(1+L)^kf\|_2$ for any $f\in H^{\infty}_L(G)$, an hence in particular for any $f\in L^2(G)$ compactly supported. We now can deduce property RD, because for $f\in L^2(G)$ supported in $B_L(R)$ for $R\geq 1$, $\|(1+L)^kf\|_2\leq C(1+R)^k\|f\|_2$.\end{proof}
%%%%%%%%%%%%%%%%%%%%%%%%%%%%%%%%%%%%%%%%%%%%%%%%%%%%%
\section{Elementary stability results of property RD}\label{stabilitefacile}
%%%%%%%%%%%%%%%%%%%%%%%%%%%%%%%%%%%%%%%%%%%%%%%%%%%%%%
Property RD is not stable under general extensions. Indeed, all abelian groups have property RD, but there are solvable groups with exponential volume growth, and these groups do not have property RD (see Proposition \ref{moyennable} below). However property RD is stable under direct products and as will be proved later under some central extensions (see Proposition \ref{RDextensionsPol}).
\begin{lem}\label{produits}Let $G_1$ and $G_2$ be two compactly generated groups, endowed with length functions $L_1$ and $L_2$. Let $G=G_1\times G_2$ and $L=L_1+L_2$. Then $(G,L)$ has property RD if and only if $(G_1,L_1)$ and $(G_2,L_2)$ do.\end{lem}
\begin{proof} For $f\in L^2(G)$ compactly supported, define
$$f_1(x)=\sqrt{\int_{G_2}|f(x,y)|^2\vec{d}y}\in L^2(G_1)$$
Then $\|f\|_{L^2(G)}=\sqrt{\int_{G_1}|f_1(x)|^2\vec{d}x}=\|f_1\|_{L^2(G_1)}$ and
$$L_F(g)=\int_Gf(y)g(y^{-1}x)\vec{d}y=\int_{G_1\times G_2}f(y_1,y_2)g(y_1^{-1}x_1,y_2^{-1}x_2)\vec{d}y_1\vec{d}y_2.$$
Now assume property RD for $G_1$ and $G_2$ with respective constants $C_1,D_1$ and $C_2,D_2$, and take $f\in L^2(G)$ compactly supported in the ball of radius $R>1$ for the length $L=L_1+L_2$. Fixing the variable $x_1$ we see that
\begin{eqnarray*}& & \int_{G_2}\left|\int_{G_1\times G_2}f(y_1,y_2)g(y_1^{-1}x_1,y_2^{-1}x_2)\vec{d}y_1\vec{d}y_2\right|^2\vec{d}x_2\\
&\leq&\left(\int_{G_1}\left(\int_{G_2}\left|\int_{G_2}f(y_1,y_2)g(y_1^{-1}x_1,y_2^{-1}x_2)dy_2\right|^2\vec{d}x_2\right)^{1/2}\vec{d}y_1\right)^2\\
&\leq & C_2^2R^{2D_2}\left|\int_{G_1}f_1(y_1)g_1(y_1^{-1}x_1)\vec{d}y_1\right|^2\end{eqnarray*}
where the first inequality is Minkowsky (see \cite{Rudin} Theorem 3.29) and the last inequality is the assumption that $(G_2,L_2)$ has property RD. Since $f_1$ is compactly supported in $B_{L_1}(R)$, integrating with respect to $x_1$ yields
$$\|L_F(g)\|_{L^2(G)}\leq C_2R^{D_2}\|L_{F_1}(g_1)\|_{L^2(G_1)}\leq CR^D\|f\|_{L^2(G)}\|g\|_{L^2(G)}$$
where $C=C_1C_2$ and $D=D_1+D_2$. The last inequality follows from the assumption that $(G_1,L_1)$ has property RD.

\medskip

Conversely, assume that $G$ has property RD and take $f_1\in{\mathcal C}_0(G_1)$, supported in $B_{L_1}(R)$ for $R\geq 1$. We fix $U$ a compact neighborhood of $1\in G_2$ contained in $B_{L_2}(1)$. Consider the function ${f}\in L^2(G)$ defined by
$${f}(y_1,y_2)=f_1(y_1){\bf 1}_U(y_2),$$
for all $(y_1,y_2)\in G_1\times G_2=G$ and where ${\bf 1}_U$ denotes the characteristic function of $U$. Let $L_U$ denote the left convolution by ${\bf 1}_U$. Since $U$ is compact, we have that $\|L_U\|_{2\ra 2}<\infty$. Notice that $\|L_F\|_{2\ra 2}\geq\|L_U\|_{2\ra 2}\|L_{F_1}\|_{2\ra 2}$, and hence we have that
\begin{eqnarray*}\|L_{F_1}\|_{2\ra 2}&\leq &\|L_U\|_{2\ra 2}^{-1}\|L_F\|_{2\ra 2}\leq \|L_U\|_{2\ra 2}^{-1}CR^D\|f\|_{L^2(G)}\\
& = &\sqrt{\vec{\nu}(U)}\|L_U\|_{2\ra 2}^{-1}CR^D\|f_1\|_{L^2(G_1)},\end{eqnarray*}
where $C$ and $D$ are the constants of the RD inequality. We conclude by Lemma \ref{EquivRD}.\end{proof}
We shall see in the sequel that property RD is not closed under passing to general subgroups but transfers to open subgroups. More precisely we have the following.
\begin{lem}\label{ssgroupeouvert} Let $(G,L)$ have property RD and take $H<G$ an open subgroup. Then $(H,L')$ has property RD, where $L'$ is the length function $L$ restricted to $H$.\end{lem}
\begin{proof}Since $H$ is open, the Haar measure on $H$ is the restriction of the one on $G$. Let $f\in L^2(H)$ be supported on $B_{L'}(R)$ for some $R\geq 1$. Extend $f$ to $\widetilde{f}\in L^2(G)$ by setting $\widetilde{f}=0$ on $G\setminus H$, so that $\|f\|_{L^2(H)}=\|\widetilde{f}\|_{L^2(G)}$ and $\widetilde{f}$ is supported on $B_L(R)$. Then
$$\|L_F\|_{2\ra 2}\leq\|L_{\widetilde{F}}\|_{2\ra 2}\leq CR^D\|\widetilde{f}\|_{L^2(G)}=CR^D\|f\|_{L^2(H)}.$$\end{proof}
The following will allow us to extend our results on connected compactly generated groups to virtually connected compactly generated groups.
\begin{lem}\label{extfin} Let $G$ be a compactly generated group, and $H$ a closed finite index subgroup. Then $G$ has property RD if and only if $H$ does.\end{lem}
\begin{proof}Since $H$ has finite index in $G$, it is an open subgroup. Hence if $G$ has property RD then so does $H$ by Lemma \ref{ssgroupeouvert}. Conversely, we start by defining the map $\bar{\ }:L^2(H)\to L^2(G)$, $f\mapsto\bar{f}$ as the extension by 0 and by noticing that for $f\in{\mathcal C}_0(H)$, then 
\begin{equation}\label{NormeOp=}\|L_{\bar{F}}\|_{2\ra 2}=\|L_F\|_{2\ra 2}.\end{equation} 
Indeed, choose a system of representatives $R$ for right $H$-cosets, and for $u\in{\mathcal C}_0(G)$, $g\in G$, write $u(g)=\sum_{r\in R}\bar{u}_r(r^{-1}g)$ where $u_r\in{\mathcal C}_0(H)$ is such that ${u}_r(h)=u(rh)$. Since $H$ is open in $G$, a left Haar measure on $G$ restricted to $H$ coincides with a left Haar measure on $H$. Hence $\|u\|_{L^2(G)}^2=\sum_{r\in R}\|u_r\|_{L^2(H)}^2$. For $x=rh\in G$ we have that
\begin{eqnarray*}R_{\bar{F}}u(x)&=&\int_Gm^{-1}(y)u(xy^{-1})\bar{f}(y)\vec{d}y\\
&=&\int_Hm^{-1}(y)u_r(hy^{-1})f(y)\vec{d}y=R_Fu_r(h),\end{eqnarray*}
where the $*$ is now the convolution in $H$. Computing the $L^2$-norm we get
\begin{eqnarray*}\|R_{\bar{F}}u\|_{L^2(G)}^2&=&\sum_{r\in R}\|R_{F}u_r\|_{L^2(H)}^2\\
&\leq &\|R_F\|_{2\ra 2}\sum_{r\in R}\|u_r\|_{L^2(H)}^2=\|R_F\|_{2\ra 2}\|u\|_{L^2(G)}^2.\end{eqnarray*}
This proves that $\|R_{\bar{F}}\|_{2\ra 2}=\|R_F\|_{2\ra 2}$ and (\ref{NormeOp=}) follows using (\ref{def}). Now assume that $H$ has property RD, and take $f\in L^2(G)$ supported on a ball of radius $R$ for some compact generating set $K\subseteq G$, where we assume that $K$ contains a finite system of representatives $x_1,\dots,x_n$ for right $H$-cosets. Since $G=\coprod_{i=1}^nx_iH$, we can write $f={f}_1+\dots+{f}_n$, where ${f}_i=f{\bf 1}_{[x_i]}$ and ${\bf 1}_{[x_i]}$ is the characteristic function of the coset $x_iH$ (the ${\bf 1}_{[x_i]}$'s are continuous since the $x_iH$'s are open and closed in $G$). Moreover, for each $i=1,\dots,n$, we can define $\tilde{f}_i\in L^2(H)$ by $\tilde{f}_i(h)=f_i(x_ih)$ and so $\|\tilde{f}_i\|_{L^2(H)}=\|f_i\|_{L^2(G)}$. Since $f_i$ is a left translation of $\bar{\tilde{f}}_i$, we have that $\|L_{\bar{\tilde{F}_i}}\|_{2\ra 2}=\|L_{F_i}\|_{2\ra 2}$ and hence $\|L_{F_i}\|_{2\ra 2}=\|L_{\tilde{F}_i}\|_{2\ra 2}$ because of the discussion in the beginning of the proof. We now compute
\begin{eqnarray*}\|L_F\|_{2\ra 2}&\leq&\sum_{i=1}^n\|L_{F_i}\|_{2\ra 2}=\sum_{i=1}^n\|L_{\tilde{F}_i}\|_{2\ra 2}\leq CR^D\sum_{i=1}^n\|\tilde{f}_i\|_{L^2(H)}\\
&\leq& CR^D\sqrt{n}\sqrt{\sum_{i=1}^n\|f_i\|^2_{L^2(G)}}=CR^D\sqrt{n}\|f\|_{L^2(G)},\end{eqnarray*}
where the first inequality is the triangle inequality, the second one is property RD and the last one is Cauchy-Schwartz.\end{proof}
{\it Remarks.} 1. Lemma \ref{extfin} does not extend to the case where $H$ is a closed subgroup of $G$ and $G/H$ is compact. The semisimple groups with finite center studied in Section \ref{ss} provide counter-examples (such a group $\Sigma=NAK$ has property RD by Theorem \ref{principal}, but $NA$ does not according to Proposition \ref{moyennable}). 

\smallskip

2. According to Jolissaint \cite{Jol} Proposition A.3, if a compactly generated group $G$ has a discrete cocompact subgroup with property RD, then $G$ has property RD as well. For instance, Jolissaint proved that discrete groups acting properly and cocompactly on Riemannian manifolds with pinched negative sectional curvature have property RD. He deduces (Corollary A.4) that $SL_2(\R)$ as well as any connected non compact Lie group of real rank one and finite center, has property RD.

\smallskip

3. It is not known whether property RD passes to cocompact lattices and this appears to be a hard question. So far only a few cocompact lattices in semisimple Lie groups are known to have property RD (see \cite{RRS}, \cite{L1} and \cite{moi}), and the methods used to establish property RD for those groups are quite different from what we do here for connected groups. It is a conjecture of Valette \cite{Val} that cocompact lattices in real and $p$-adic semisimple Lie groups have property RD. 
%%%%%%%%%%%%%%%%%%%%%%%%%%
\section{Amenability}\label{amen}
%%%%%%%%%%%%%%%%%%%%%%%%%
One of the generalizations of Kesten's theorem concerning amenability reads as follows (see \cite{L}). In the notation introduced above, $G$ is amenable if and only if for any $\mu$ a (non-negative) Borel measure on $G$
\begin{equation}\label{KestenL}\|L_\mu\|_{2\ra 2}=\mu(G)=\|\mu\|.\end{equation} 
In particular, if $G$ is amenable and $\mu$ is a (non-negative) Borel measure such that $\int_G m(y)^{-1/2}d\mu(y)<\infty$ we have
\begin{equation}\label{KestenR}\|R_{\mu}\|_{2\ra 2}=\int_G m(y)^{-1/2}d\mu(y)<\infty.\end{equation}
We say that $(G,L)$ is of polynomial volume growth $d$ if there exists a constant $c$ such that $\vec{\nu}(B_L(n))\le cn^d$. If $G$ is not unimodular, i.e., $m\not\equiv 1$ then $\vec{\nu}(B_L(n))$ grows exponentially (see e.g. \cite{VSC} Chapter IX).
\begin{defin}
Let $L$ be a length function on a locally compact group $G$. A function $f$ is \emph{radial} (with respect to $L$) if $L(x)=L(y)$ implies $f(x)=f(y)$. We say that $(G,L)$ has \emph{property RD-rad} (respectively \emph{property RD'-rad}) if $(G,L)$ has property RD$_E$ (respectively property RD'$_E$) for $E$ the set of radial functions (with respect to $L$) in $L^2(G)$.\end{defin}
The following is a straightforward adaptation of Jolissaint's Corollary 3.1.8. of \cite{Jol}.
\begin{pro}\label{moyennable}Let $G$ be a locally compact, amenable group and $L$ a proper length function on $G$. Then the following are equivalent:
\begin{itemize}\item[(1)]$(G,L)$ has polynomial volume growth.
\item[(2)]$(G,L)$ has property RD.
\item[(3)]$(G,L)$ has property RD'-rad.\end{itemize}
\end{pro}
\begin{proof}As $G$ is amenable, we have
\begin{equation}\label{CaSc}\|L_F\|_{2\ra 2}=\int_G f(x)\vec{d}x\le\sqrt{\vec{\nu}(\mbox{supp}(F))}\|f\|_2,\end{equation}
the first equality follows from equation (\ref{KestenL}) and the inequality by Cauchy-Schwartz. We prove that (1) implies (2). If $G$ has polynomial volume growth of degree $d$ and $f$ is supported on $B_L(R)$ for some $R\geq 1$, then 
$$\sqrt{\vec{\nu}(\mbox{supp}(F))}\le \vec{\nu}(B_L(R))\le cR^d$$ 
and property RD' (hence RD) follows. As RD and RD' are equivalent (see Corollary \ref{RD<=>RD'}), (2) obviously implies (3). We now show that (3) implies (1). Let $R\geq 1$, applying property RD'-rad to the radial function $f={\mathbf 1}_{B_L(R)}$ and using the equality in (\ref{CaSc}) gives
$$\vec{\nu}(B_L(R))\le CR^D\sqrt{\vec{\nu}(B_L(R))},$$
where $C,D$ are the constants of Definition \ref{DefRDRD'}. It follows that 
$$\vec{\nu}(B_L(R))\leq C^2R^{2D}.$$
\end{proof}
\noindent
{\it Remark.} We shall see in Section \ref{NAaRDrad} below that RD-rad does not imply polynomial volume growth. In particular, property RD-rad is in general not equivalent to property RD'-rad.
%%%%%%%%%%%%%%%%%%%%%%%%%%%%%%%%%%%%%%%%%%%%%%%%%%%%%%%
\section{Property RD-rad on Iwasawa $NA$ groups}\label{NAaRDrad}
%%%%%%%%%%%%%%%%%%%%%%%%%%%%%%%%%%%%%%%%%%%%%%%%%%%%%%%%%%
In this section we will see that amenable groups with exponential volume growth may satisfy property RD-rad. More precisely we prove the following.
\begin{theo}\label{NAaRDradTh} Let $\Sigma=NAK$ be a connected non-compact semisimple Lie group with finite center. Let $d$ be the canonical Riemannian distance on the symmetric space $X=\Sigma/K$. Let $S$ be the solvable group $S=NA$ identified as a manifold with $X$ and let $L(g)=d(K,gK)$ for $g\in S$. Then $S$ has property RD-rad with respect to $L$.
\end{theo}
The proof of this theorem relies on results in \cite{CGHM} that we now explain. Let us start by recalling some standard facts concerning semisimple Lie groups (we follow the notations of \cite{CGHM}). Let ${\bfsigma}$ be the Lie algebra of $\Sigma$ and $\theta$ be a Cartan involution with associated Cartan decomposition ${\bfsigma}={\mathbf k}\oplus{\mathbf p}$. Fix a maximal abelian subalgebra ${\mathbf a}$ in ${\mathbf p}$ (see Chapter IX of \cite{Hel}). This determines a root space decomposition and, after ordering of the roots, an Iwasawa decomposition ${\bfsigma}={\mathbf n}\oplus{\mathbf a}\oplus{\mathbf k}$. Let $d$ be the canonical $\Sigma$-invariant Riemannian distance on the symmetric space $X=\Sigma/K$ (the projection of the left-invariant and $K$-right invariant canonical distance on $\Sigma$). Let $S$ be the solvable group $S=NA$. The action of the group $S$ on $X$ is transitive with trivial stabilizers so that we can identify $X$ and $S$ as manifolds and under this identification the $\Sigma$-invariant measure on $X$ is the left Haar measure on $S$. We can define the ``modular'' function $m$ on $X$ by setting $m(x)=m(g)$ for the unique $g\in S$ such that $go=x$. The canonical Riemannian distance on $X$ corresponds to a left invariant Riemannian distance on $S$, which gives the length function $L$ of Theorem \ref{NAaRDradTh}. The crucial point is the following result from \cite{CGHM}.
\begin{lem}\label{crucial} Let $E$ be a $K$-invariant measurable subset of $X$ and $f$ a function in $L^2(X)$ such that $m^{1/2}f$ is $K$-invariant. Then
\begin{equation}\|{\mathbf 1}_E f\|_2=\|{\mathbf 1}_E m^{1/2} f\|_2\end{equation}
and 
\begin{equation}\|{\mathbf 1}_E f\|_1\le \|{\mathbf 1}_E\phi_0\|_2\|{\mathbf 1}_E f\|_2\end{equation}
where $\phi_0$ is the basic elementary spherical function.\end{lem}
We will need the Cartan decomposition in the form of the integrals equality
\begin{equation}\label{integrationformula}\int_S f(x)\vec{d}x=C\int_{{\mathbf a}^+}\int_K f(k\exp(H))dk D(H)dH\end{equation}
(for all $f\in L^1(G)$) where $dk$ is the normalized Haar measure on $K$, ${\mathbf a}^+$ is the positive Weyl chamber in ${\mathbf a}$ and $D(H)$ has an explicit expression in terms of the roots and satisfy (see \cite{CGHM}):
\begin{equation}\label{D(H)estimate}D(H)\le C_1\left(\frac{|H|}{1+|H|}\right)^{n-\ell}e^{2\rho (H)}.\end{equation}
Here $|H|$ is the length of $H$ in the Lie algebra ${\mathbf a}$ equipped with its canonical Euclidean structure given by the Killing form; $\rho(H)$ is the usual half sum of the positive roots counting multiplicity, $n$ is the topological dimension of $S=NA$ and $\ell$ is the dimension of $A$, i.e., the real rank of $\Sigma$. The modular function on $S=NA$ is given by $m(n\exp H)=\exp(-2\rho(H))$ in standard notation. We shall need the following estimate from \cite{CGHM}.
\begin{lem}\label{estimeRad}For $r>1$,
\begin{equation}\label{SpheriquePol}\int_{B(r)}|\phi_0(x)|^2\vec{d}x\le C r^\gamma\end{equation}
where $\gamma=2b+\ell=2\times\sharp\{\hbox{indivisible positive roots}\}+\dim(A)$.\end{lem}
\begin{proof}This follows from the above integration formula (\ref{integrationformula}), the estimate (\ref{D(H)estimate}) on $D(H)$, the fact that $\phi_0(\exp(H))\le C(1+|H|)^b e^{-\rho(H)}$, as well as $d(K,\exp(H)K)=|H|$.\end{proof}
We can now proceed with the proof of Theorem \ref{NAaRDradTh}
\begin{proof}[Proof of Theorem \ref{NAaRDradTh}]Let $B(r)=B_L(r)$ denote the ball of radius $r$ around the origin $e=e_S=e_\Sigma K$ in $S=X$. The distance $d=d_X$ on $X$ is $K$-invariant and hence, as $K$ fixes the origin in $X$, the set $B(r)$ is $K$-invariant. Now, fix a continuous function $f$ on $S$ which is radial and supported on the ball $B(r)$. Set $dF(x)=f(x)\vec{d}x$. By the inequality (\ref{def}) we  have
$$\|R_F\|_{2\ra 2}\le\int_S m(x)^{-1/2} f(x)\vec{d}x=\int_X {\mathbf 1}_{B(r)}(x)m(x)^{-1/2} f(x)\vec{d}x.$$
As $B(r)$ is $K$-invariant and $f=m^{1/2}m^{-1/2}f$ is radial, the crucial Lemma \ref{crucial} gives
$$\int_X {\mathbf 1}_{B(r)}(x)m(x)^{-1/2} f(x)\vec{d}x\le \left(\int_{B(r)} |\phi_0(x)|^2\vec{d}x\right)^{1/2}\left(\int_{B(r)} |f(x)|^2\vec{d}x\right)^{1/2}.$$
Finally, combining the above inequality with the estimate (\ref{SpheriquePol}) we get
$$\|R_F\|_{2\ra 2}\le Cr^{\gamma/2} \|f\|_2$$
as desired.\end{proof}
\noindent
{\it Remark.} Estimates similar to the ones above and implying property RD-rad on Damek-Ricci $NA$ groups have been obtained by Mustapha in \cite{M}, equation (3.8). 
%%%%%%%%%%%%%%%%%%%%%%%%%%%
\section{Semisimple groups with finite center}\label{ss}
%%%%%%%%%%%%%%%%%%%%%%%%%%%
The aim of this section is to prove the following.
\begin{theo}\label{principal} Connected non-compact semisimple Lie groups with finite center have property RD with respect to the length function associated to the canonical left-invariant Riemannian metric. In particular, they have property RD.\end{theo}
The idea is to reduce this result to Theorem \ref{NAaRDradTh}. To do so, we start with the following easy lemma (this is a classical observation used in particular in the study of Kunze-Stein phenomena, see \cite{Cow}).
\begin{lem}\label{Kinvariance} Let $G$ be a locally compact group and $K$ be a compact subgroup. Let $f\in{\mathcal C}_0(G)$ and $dF=fdx$. We set
\begin{eqnarray*}f_K(x)=\left(\int_K |f(xk)|^2 dk\right)^{1/2},\widetilde{f}(x)=\left(\int_K \int_K |f(kxk')|^2 dkdk'\right)^{1/2},\end{eqnarray*}
where $dk$ denotes the normalized Haar measure on $K$. Then $f_K,\widetilde{f}\in L^2(G)$, $\|f_K\|_2=\|\widetilde{f}\|_2=\|f\|_2$ and
\begin{eqnarray*}\|L_F\|_{2\ra 2}& \le &\|L_{F_K}\|_{2\ra 2}\\
\|L_F\|_{2\ra 2}&\le & \|L_{\widetilde{F}}\|_{2\ra 2},\end{eqnarray*}
where $dF_K=f_K\vec{d}x$ and $d\widetilde{F}=\widetilde{f}\vec{d}x$.\end{lem}
\begin{proof} The equality of the norms of $f_K,\widetilde{f}$ and $f$ follows from the fact that $m(k)=1$ for $k\in K$ (indeed, $m(k^n)=m(k)^n$ is bounded and away from 0 for all $n\in\N$ since $K$ is compact). For any function $g\in{\mathcal C}_0(G)$ and any $x\in G,k\in K$ we have that
$$F*g(x)=\int_G f(xy)g(y^{-1})\vec{d}y=\int_G f(xyk)g(k^{-1}y^{-1})\vec{d}y.$$ 
Hence,
\begin{eqnarray*}
\left|F*g(x)\right|& = & \left|\int_K\int_G f(xyk)g(k^{-1}y^{-1})\vec{d}ydk\right|\\
&\le & \int_G \left(\int_K|f(xyk)|^2dk\right)^{1/2}\left(\int_K|g(ky^{-1})|^2dk\right)^{1/2}\vec{d}y\\
& = & \int_G f_K(xy)g_K(y^{-1})\vec{d}y=F_K*g_K(x)\end{eqnarray*}
It follows that $\|L_F\|_{2\ra 2}\le\|L_{F_K}\|_{2\ra 2}$. Now we set
$$\!\!\phantom{f}_K\! f(x)=\left(\int_K |f(kx)|^2 dk\right)^{1/2},$$
and check that $(\check{F}_K)\!\!{\check{\phantom{f}}}=\!\!\!\phantom{f}_K\! F$. Applying the above inequality twice, we obtain
\begin{eqnarray*}\|L_F\|_{2\ra 2}&=&\|L_{\check{F}}\|_{2\ra 2}\leq\|L_{\check{F}_K}\|_{2\ra 2}=\|L_{(\check{F}_K)\!\!{\check{\phantom{f}}}}\|_{2\ra 2}\\
&=&\|L_{\!\!\!\phantom{f}_K\! F}\|_{2\ra 2}\leq\|L_{\!\!\!\phantom{f}_K\! F_K}\|_{2\ra 2}=\|L_{\widetilde{F}}\|_{2\ra 2}.\end{eqnarray*}
\end{proof}
{\it Remark.} In Lemma \ref{Kinvariance} we also have $\|R_F\|_{2\ra 2}\le\|R_{F_K}\|_{2\ra 2}$ and $\|R_F\|_{2\ra 2}\le\|R_{\widetilde{F}}\|_{2\ra 2}$. Indeed, $\|R_{F}\|_{2\ra 2}=\|L_{\hat{F}}\|_{2\ra 2}$ for any finite measure $F$ with $dF=f\vec{d}x$ by (\ref{def}), moreover, $(F_K)\!{\hat{\phantom{f}}}=\hat{F}_K$ and $\widetilde{\hat{F}}=\hat{\widetilde{F}}$.

\medskip

Let us now recall some facts about homogeneous spaces, see, e.g., \cite{WS2,WS3}. Let $G$ be a locally compact group which acts continuously, transitively on a space $X$ with compact stabilizers. Fix $o\in X$ and let $K$ denote the stabilizer of $o$ so that $X=G/K$. For $x\in X$, let $\bar{x}$ be an element of $G$ such that $\bar{x}o=x$. Let $p(x,y)$ be a locally integrable non-negative kernel which is $G$-invariant (i.e. $p(gx,gy)=p(x,y)$ for any $g\in G$) and such that the support of $p(x,\cdot)$ is compact for all $x\in X$. Let $dx$ denote the $G$-invariant measure on $X$ so that $\vec{d}g=dxdk$ where $dk$ is the normalized Haar measure on $K$. Set 
$$\phi(g)=p(go,o)=p(o,g^{-1}o)\hbox{ and }d\Phi(g)=\phi(g)\vec{d}g.$$
Note that $\phi$ satisfies $\phi(gk)=\phi(kg)=\phi(g)$ for all $g\in G$, $k\in K$. Moreover one checks that the right convolution operator $R_\Phi$ ``realizes'' on $G$ the operator $P:{\mathcal C}(X)\to{\mathcal C}(X),f\mapsto Pf$ where
$$Pf(x)=\int_X p(x,y)f(y)dy.$$
More precisely, define 
\begin{eqnarray*}S:{\mathcal C}(G)&\to & {\mathcal C}(X)\\
 f&\mapsto & Sf:\{x\mapsto\int_K f(\bar{x}k)dk\}\\
T:{\mathcal C}(X)&\to & {\mathcal C}(G)\\
f&\mapsto & Tf:\{g\mapsto f(go)\}.\end{eqnarray*}
Then $S$ and $T$ act on the respective $L^2$-spaces without increasing norms, and $P=SR_\Phi T$, $R_\Phi=TPS$. In particular, $\|P\|_{2\ra 2}=\|R_\Phi\|_{2\ra 2}$ and
$$\int_G|\phi(g)|^2\vec{d}g= \int_X |p(x,o)|^2 dx.$$
Given $G$, $X$ and $p$ as above, if $Q$ is another locally compact group which acts continuously and transitively on $X$ with compact stabilizers and such that $p$ is $Q$-invariant, we get right convolution operators $R^G_{\Phi^G}$ and $R^Q_{\Phi^Q}$ on $G$ and $Q$ respectively with
\begin{equation}\label{sousgroupe}\|P\|_{2\ra 2}=\|R^G_{\Phi^G}\|_{2\ra 2}=\|R^Q_{\Phi^Q}\|_{2\ra 2}\hbox{ and }\|\Phi^G\|_2=\|\Phi^Q\|_2.\end{equation}
The following will be used later to extend our results on connected Lie groups to compactly generated groups.
\begin{lem}\label{ExtCompacte}Let $1\to H\to G\to Q\to 1$ be a short exact sequence of compactly generated groups, and assume that $H$ is compact. Then $G$ has property RD if and only if $Q$ has property RD.\end{lem}
\begin{proof}First assume that $G$ has property RD and take $f\in L^2(Q)$ supported on $B(R)$ for $R\geq 1$ and a compact symmetric neighborhood of the identity $A$ generating $Q$. Set ${f}^{\sharp}=f\circ\pi$, where $\pi:G\to Q$ is the projection with kernel $H$ in the short exact sequence above. Then ${f}^{\sharp}$ has its support in $B(R+1)$ for the compact generating set $\pi^{-1}(A)$. Choosing the normalized Haar measure on $H$ we get, using equality (\ref{sousgroupe}) above
$$\|R_F\|_{2\ra 2}=\|R_{{F}^{\sharp}}\|_{2\ra 2}\leq C(R+1)^D\|{f}^{\sharp}\|_{L^2(G)}=C(R+1)^D\|f\|_{L^2(Q)}.$$
Conversely, assume that $Q$ has property RD and take $f\in L^2(G)$ supported on $B(R)\subseteq G$ for $R\geq 1$ and a compact generating set as above. Define $f_H\in L^2(Q)$ as in Lemma \ref{Kinvariance}, so that $f_H$ is supported on $B(R)\subseteq Q$, and
$$\|R_F\|_{2\ra 2}\leq\|R_{F_H}\|_{2\ra 2}\leq CR^D \|f_H\|_{L^2(Q)}=CR^D\|f\|_{L^2(G)}.$$\end{proof}
\begin{proof}[Proof of Theorem \ref{principal}]As in Theorem \ref{NAaRDradTh}, let $\Sigma=NAK$ be a connected non-compact semisimple Lie group with finite center. We will show that $\Sigma$ equipped with its canonical $K$-bi-invariant Riemannian metric has property RD. According to Lemma \ref{Kinvariance} we can consider $K$-bi-invariant functions only. For $K$-bi-invariant functions, equality (\ref{sousgroupe}) reduces the situation to $K$-invariant functions on $X$, or equivalently to radial functions on $S=NA$. More precisely, there is a constant $C$ such that for $f$ with compact support in the ball of radius $R$ and $\widetilde{f}$ as in Lemma \ref{Kinvariance} we have
$$\|L_F\|_{2\ra 2}\le\|L_{\widetilde{F}}\|_{2\ra 2}=\|R_{\widetilde{F}}\|_{2\ra 2}\le CR^{\gamma/2}\|f\|_2,$$
where $\gamma$ is as in Lemma \ref{estimeRad}. The first inequality holds because of Lemma \ref{Kinvariance}, the middle equality because $\Sigma$ is unimodular and the last inequality holds because of Theorem \ref{NAaRDradTh} combined with equality (\ref{sousgroupe}).\end{proof}
%%%%%%%%%%%%%%%%%%%%%%%%%%%%%%%%%%%%%%%%%%%%%%%%%
\section{Central extensions and property RD}\label{ext}
%%%%%%%%%%%%%%%%%%%%%%%%%%%%%%%%%%%%%%%%%%%%%%%%%%%
The aim of this section is to investigate the stability of property RD under certain central extensions. Recall that if $E$ is a central extension of a group $G$ by a group $A$
$$1\to A\to E\to G\to 1,$$
then $E$ is $A\times G$ as a set, and the group law is given by
$$(a,g)(a',g')=(a+a'+c(g,g'),gg')$$
where $c:G\times G\to A$ is such that $c(g,1)=c(1,g)=0\in A$ (we denote by ``$+$'' the group law in $A$), and for $g,g',g''$ the cocycle relation is satisfied:
$$c(g,g')+c(gg',g'')=c(g',g'')+c(g,g'g'').$$
The cocycle $c$ satisfies
$$(c(g,g'),1)=\sigma(g)\sigma(g')\sigma(gg')^{-1},$$
where $\sigma:G\to E$ is the section $g\mapsto (0,g)$ of $p:E\to G, (a,g)\mapsto g$.
\begin{defin}Let $A$ and $G$ be two compactly generated groups, with $A$ abelian. Let $\Gr:{\N}\to{\N}$ be a function and let $c:G\times G\to A$ be a 2-cocycle. We say that $c$ has \emph{growth at most $\Gr$} if 
$$\sup\{L_A(c(g,g'))|g,g'\in B(n)\subseteq G\}\leq\Gr(n).$$
We say that a 2-cocycle has \emph{polynomial growth} if $\Gr$ is a polynomial.\end{defin}
The following shows that a central extension defined by a cocycle with polynomial growth of a group with property RD has again property RD.
\begin{pro}\label{RDextensionsPol}Let $1\to A\to E\to G\to 1$ be an exact sequence of compactly generated groups with $A$ closed and central. Assume there exists a compact neighborhood $U=U^{-1}$ of 1 in $G$ which generates $G$ and a Borel section $\sigma:G\to E$ of the canonical projection $p:E\to G$ such that $\sigma(1)=1$ and such that $\sigma(U)$ is compact and such that the cocycle defined by $c(g,g')=\sigma(g)\sigma(g')(\sigma(gg'))^{-1}$ has polynomial growth. If $G$ has property RD then $E$ has property RD as well.\end{pro}
The map $A\times G\to E$, $(a,g)\mapsto a\sigma(g)$ is a Borel isomorphism because $\sigma$ is Borel. We identify $E$ with $A\times G$ and $\sigma(g)$ with $(1,g)$ as explained above. Assume that the defining cocycle $c$ has polynomial growth $\Gr(n)$. Let $U=U^{-1}$ and $S=S^{-1}$ be compact generating sets for $G$ and $A$ respectively. By hypothesis we may choose $U$ such that $\sigma(U)$ is compact in $E$. Then 
$$T=\{(s,u)|s\in S,u\in U\}\cup\{(s,u)|s\in S,u\in U\}^{-1}$$ 
is a compact generating set for $E$. To see that it is a neighborhood of 1, let $K$ be a compact neighborhood of 1 in $E$. If we choose $S\subseteq A$ big enough so that $K\sigma(U)^{-1}\cap A\subseteq S$, then $T$ as above is a compact neighborhood of 1 because $K\cap p^{-1}(U)=K\cap A\sigma(U)\subseteq S\sigma(U)\subseteq T$. One checks that $T$ generates $E$. We start with the following observation.
\begin{lem}\label{dehors}Under the above assumptions, there are constants $C$ and $k$ such that for any $a\in A$ then $L_S(a)\leq C(1+L_T(a,1))^{k}$.\end{lem}
\begin{proof}Take $(a,1)$ of length $r$ in $E$ and write $(a,1)=(a_1,x_1)\dots (a_r,x_r)$ where $(a_i,x_i)$ belong to the generating set $T$, for $i=1,\dots,r$. Note that $L_U(x_i)=1$ but $L_S(a_i)\leq M=\sup_{u\in U}L_S(c(u,u^{-1}))+1$. Then, defining $\lambda_i=x_1\dots x_i$ for $i=1,\dots,r$ we get:
$$(a,1)=(a_1,x_1)\dots (a_r,x_r)=(\sum_{i=1}^ra_i+\sum_{i=1}^{r-1}c(\lambda_i,x_{i+1}),1),$$
so that
$$L_S(a)\leq\sum_{i=1}^rL_S(a_i)+\sum_{i=1}^{r-1}L_S(c(\lambda_i,x_{i+1}))\leq  rM + (r-1)\Gr(r)\leq C(1+r)^{k}$$
choosing $C$ and $k$ appropriately.\end{proof}
\begin{proof}[Proof of Proposition \ref{RDextensionsPol}]As $G$ has property RD, it is unimodular, and it follows from the Corollaire in Chapitre VII, Paragraphe 2, Num\'ero 7 of \cite{BouIntegration} that $E$ is also unimodular. First notice that a compactly generated abelian group is automatically of polynomial growth for any algebraic length, and thus has property RD. Denote by $C_G, D_G$ and $C_A, D_A$ the constants needed for the RD inequality (as in Definition \ref{DefRDRD'}) for $G$ and $A$ respectively. For $f,g\in {\mathcal C}_0(E)$, we define $f_y(a)=f(a,y)$ and $g'_{(y,x)}(a)=g_{y^{-1}x}(a-c(y,y^{-1}x))$. For almost all $x,y$, the elements $f_y$ and $g'_{(y,x)}$ belong to $L^2(A)$ and since $c$ is measurable we have that
\begin{eqnarray*}f*g(a,x)=\int_G\left(\int_Af_y(b)g'_{(y,x)}(a-b)db\right)dy=\int_Gf_y*g'_{(y,x)}(a)dy.\end{eqnarray*}
Now squaring and integrating over $E$ the above expression gives
\begin{eqnarray*}\|f*g\|^2_{L^2(E)}&=&\int_G\left(\int_A\left|\int_Gf_y*g'_{(y,x)}(a)dy\right|^2da\right)^{\frac{1}{2}2}dx\\
&\leq & \int_{G}\left(\int_{G}\|f_y*g'_{(y,x)}\|_{L^2(A)}dy\right)^2dx.\end{eqnarray*}
Now assume that the support of $f$ is contained in the ball of radius $r$ and for $y\in{G}$, let us look at the support of $f_y$. Take $a$ in the support of $f_y$, then $L_T(a,y)\leq r$, so that $L_T(a,1)\leq L_T(a,y)+L_T(1,y)\leq 2r$ and thus Lemma \ref{dehors} implies that $L_S(a)\leq C(1+2r)^{k}$. Now we can proceed with the computation applying property RD for $A$ to $\|f_y*g'_{(y,x)}\|_{L^2(A)}$ and so
$$\|f*g\|^2_{L^2(E)}\leq \int_{G}\left(\int_{G}C_A(C(1+2r)^k)^{D_A}\|f_y\|_{L^2(A)}\|g'_{(y,x)}\|_{L^2(A)}dy\right)^2dx.$$
Finally, define $\tilde{f},\tilde{g}\in L^2(G)$ by $\tilde{f}(y)=\|f_y\|_{L^2(A)}$ and $\tilde{g}(y)=\|g_y\|_{L^2(A)}$, so that clearly $\|\tilde{f}\|_{L^2(G)}=\|f\|_{L^2(E)}$ and $\|\tilde{g}\|_{L^2(G)}=\|g\|_{L^2(E)}$. Notice that $\tilde{f}$ is supported on the ball of radius $r$. Indeed, if $y$ is in the support of $\tilde{f}$, then there exists $b\in A$ such that $(b,y)$ is in the support of $f$, which is contained in the ball of radius $r$ in $E$. Writing $(b,y)=(a_1,y_1)\dots (a_r,y_r)$ we see in particular that $y=y_1\dots y_r$, i.e. the length of $y$ in ${G}$ (with respect to the generating set $U$) is shorter than $r$. Concerning $g$, we have that:
$$\|g'_{(y,x)}\|_{L^2(A)}=\|g_{y^{-1}x}\|_{L^2(A)}=\tilde{g}(y^{-1}x)$$
(we performed the change of variable $a\mapsto a-c(y,y^{-1}x)$). Going back to the computation of $\|f*g\|^2_{L^2(E)}$ we now get that:
\begin{eqnarray*}\|f*g\|^2_{L^2(E)}&\leq&C^2_{A}(C(1+2r)^{k})^{2D_A}\|\tilde{f}*\tilde{g}\|_{L^2(G)}^2\\
&\leq&C^2_{A}(C(1+2r)^{k})^{2D_A}C_{G}^2r^{2D_G}\|f\|_{L^2(E)}^2\|g\|_{L^2(E)}^2\\
&\leq & C_E^2r^{2D_E}\|f\|_{L^2(E)}^2\|g\|_{L^2(E)}^2,\end{eqnarray*}
choosing $C_E=3^{kD_A}C_AC^{D_A}C_G$ and $D_E=kD_A+D_G$. We conclude that $E$ has property RD by Lemma \ref{EquivRD} (2) together with density of ${\mathcal C}_0(E)$ in $L^2(E)$.\end{proof}
%%%%%%%%%%%%%%%%%%%%%%%%%%%%%%%%%%%%%%%%%%%%%%%
\section{Lie groups with property RD}\label{gss}
%%%%%%%%%%%%%%%%%%%%%%%%%%%%%%%%%%%%%%%%%%%%%%%
In this section we shall prove the following (i.e., the implication (c)$\Rightarrow$(a) of our main theorem).
\begin{theo}\label{NonB=>RD} Let $G$ be a connected Lie group such that its universal cover $\widetilde{G}$ decomposes as $\widetilde{S}\times\widetilde{Q}$, where $\widetilde{S}$ is semisimple and $\widetilde{Q}$ has polynomial growth. Then $G$ has property RD.\end{theo}
\begin{cor}Semisimple Lie groups have property RD.\end{cor}
The idea is to show that a group $G$ as in Theorem \ref{NonB=>RD} is a central extension defined by a Borel cocycle of polynomial growth of a group with property RD by a compactly generated abelian group. We start with the following.
\begin{pro}\label{ExtendBorelSection} Let $p:E\to G$ be a surjective homomorphism of compactly generated groups. If $p$ admits a local Borel section $\sigma_K:K\to E$ defined on a compact symmetric generating neighborhood $K$ of 1 and such that $\sigma_K(K)$ is relatively compact in $E$, then there is a Borel section $\sigma:G\to E$ which extends $\sigma_K$ and which is Lipschitz with respect to the algebraic lengths.\end{pro}
For the proof we need the following.
\begin{lem}\label{partition}Let $G$ be a compactly generated group and $K$ be a compact symmetric neighborhood of 1 generating $G$. Then there is a countable \emph{pointed partition} $(G_n,g_n)$, that is a partition
$$G=\coprod_{n\in{\N}}G_n,$$
where the $G_n$'s are Borel subsets of $G$ and $g_n\in G_n$, such that $g_n^{-1}G_n\subseteq K$.\end{lem}
\begin{proof}Let $\{g_n\}\subseteq G$ be a maximal subset of elements with the property that $d(g_n,g_m)=L_K(g_n^{-1}g_m)>1$. Notice that since the ball of radius 1 is a neighborhood of 1, the set of $g_n$'s is discrete in $G$. Since a ball of finite radius is compact, there are only finitely many $g_n$'s in each ball of finite radius, so countably many altogether. Since $\{g_n\}$ is maximal, the union of balls of radius 1 centered at the $g_n$'s cover $G$ (if not, then there would be $g\in G$ not in $\{g_n\}$ and at distance greater than 1 to any $g_n$, which contradicts maximality). Denote by $B(g_n,r)=g_nB(r)$ the ball of radius $r$ centered at $g_n$. We define the $G_n$'s as follows:
\begin{eqnarray*}G_0=K=B(1),G_1=B(g_1,1)\setminus G_0,\,\dots\ G_n=B(g_n,1)\setminus (\bigcup_{k< n}G_k),\,\dots\end{eqnarray*}
It is a partition of $G$ by construction, and $g_n\in G_n$ because for any $n\ne m$ we have that $d(g_n,g_m)>1$, so that $g_n\not\in B(g_m,1)$. Finally, $g_n^{-1}G_n\subseteq g_n^{-1}B(g_n,1)=g_n^{-1}g_nK=K$ and the proof is complete.\end{proof}
\begin{proof}[Proof of Proposition \ref{ExtendBorelSection}]Let $K$ be as in the proposition and let $(G_n,g_n)$ be a pointed partition of $G$ as in Lemma \ref{partition}. Let $S$ be a compact symmetric generating set for $E$. For each $n\in{\N}$, let $e_n\in E$ be a pre-image of $g_n$ of minimal length in the alphabet $S$. Define 
$$\sigma_n:G_n\to E, x\mapsto e_n\sigma_K(g_n^{-1}x)$$
and $\sigma:E\to E$ by $\sigma=\coprod_{n\in{\N}}\sigma_n$, so that $\sigma$ is a Borel map. We check that it is a section for $p$. For $g\in G_n$, we have
$$p\sigma(g)=p\sigma_n(g)=p(e_n\sigma_K(g_n^{-1}g))=p(e_n)p\sigma_K(g_n^{-1}g)=g_ng_n^{-1}g=g.$$
Now let us prove that the section $\sigma$ we just obtained is Lipschitz. Let $C=\sup\{L_S(g)|g\in\sigma_K(K)\}$. Since $\sigma_K(K)$ is relatively compact in $E$ we have that $C<\infty$. For $g_n$ of length $m$ if we write $g_n=k_1\dots k_m$ with all $k_i\in K$ we have that $L_E(\sigma_K(k_1)\dots\sigma_K(k_m))\leq Cm$ and $p(\sigma_K(k_1)\dots\sigma_K(k_m))=g_n$. Since $e_n$ is a shortest pre-image of $g_n$ we deduce 
$$L_S(e_n)\leq L_S(\sigma_K(k_1)\dots\sigma_K(k_m))\leq Cm=CL_K(g_n).$$ 
Finally, take $g\in G$ and $n\in{\N}$ such that $g\in G_n$. We have
\begin{eqnarray*}L_S(\sigma(g))&=&L_S(e_n\sigma_K(g_n^{-1}g))\leq L_S(e_n)+C\\
&\leq & CL_K(g_n)+C\leq C(L_K(g)+2),\end{eqnarray*}
since $L_K(g_n^{-1}g)=L_K(g^{-1}g_n)\leq 1$ because $g_n^{-1}g\in K$ if $g\in G_n$.\end{proof}
\begin{defin}Let $D:{\N}\to{\N}$ be a non-decreasing function, $A<E$ be two compactly generated groups. We say that $A$ has \emph{distortion at most $D$} if there are two compact generating sets $U$ and $S$ for $A$ and $E$ respectively such that for all $a\in A$
$$L_U(a)\leq D(L_S(a))L_S(a).$$
We say that $A$ has \emph{polynomial distortion} if $D$ can be chosen to be a polynomial, and \emph{undistorted} if $D$ can be chosen constant.\end{defin}
Notice that our definition of distortion is equivalent to Gromov's one given in Chapter 3 of \cite{GNR}, as he defines (under the hypothesis of the above definition) the distortion function as
$$\hbox{\sc Disto}(r):=\frac{{\rm diam}_A(A\cap B_E(r))}{r},$$
and one easily checks that $A$ has distortion at most {\sc Disto} because $2D(n)\geq\hbox{\sc Disto}(n/2)$. Distortion of central subgroups and growth of cocycles are related as follows.
\begin{pro}\label{CocycleEtDistorsion} Let $1\to A\to E\to G\to 1$ be a central extension of compactly generated groups, and assume that $A$ has distortion at most $D$ in $E$. If $p$ admits a local Borel section $\sigma_K:K\to E$ defined on a compact symmetric generating neighborhood $K$ of 1 such that $\sigma_K(K)$ is relatively compact in $E$, then $E$ can be defined by a cocycle of growth $\Gr(r)=D(Cr)Cr$ for some constant $C$.\end{pro}
\begin{proof}Assume that $E$ is compactly generated by $S=S^{-1}$ containing 1, that $A$ is compactly generated by $U\subseteq S$ and let $K=p(S)$, it is a compact symmetric neighborhood of 1 generating $G$. We denote by $L_A,L_E$ and $L_G$ the respective length functions on $A,E$ and $G$. According to Proposition \ref{ExtendBorelSection} there is $\sigma:G\to E$ a Lipschitz section for $p:E\to G$. We choose $c(g,g')=\sigma(g)\sigma(g')\sigma(gg')$. For any $g,g'\in B(r)\subseteq G$ we have
\begin{eqnarray*}L_E(c(g,g'))&=&L_E(\sigma(g)\sigma(g')\sigma(gg')^{-1})\\
&\leq &L_E(\sigma(g))+L_E(\sigma(g'))+L_E(\sigma(gg')^{-1})\\
&\leq &C'(L_G(g)+L_G(g')+L_G(gg'))\leq Cr,\end{eqnarray*}
for some constant $C$, so that 
$$L_A(c(g,g'))\leq D(L_E(c(g,g')))L_E(c(g,g'))\leq D(Cr)Cr.$$\end{proof}
\begin{lem}\label{ChP}Let $Z$ be the center of a simply connected semisimple Lie group $\widetilde{G}$. Then $Z$ is undistorted in $\widetilde{G}$.\end{lem}
\begin{proof} Let $G=\widetilde{G}/Z$ and denote by $p:\widetilde{G}\to G$ the canonical projection. As $\widetilde{G}$ is semisimple, $Z$ is discrete and this implies that $G$ has trivial center. Let $G=NAK$ be an Iwasawa decomposition (see \cite{Hel}), since $G$ has trivial center $K$ is compact. We denote by $S$ the simply connected group $S=NA$, by $\widetilde{K}=p^{-1}(K)$ and by $\widetilde{S}$ the connected component of 1 in $p^{-1}(S)$. Consider the map 
\begin{eqnarray*}\varphi:G&\to & S\times K\\
 g&\mapsto& (s,k).\end{eqnarray*}  
On $G$ we fix a left-invariant Riemannian metric. We consider $S\times K$ as the direct product of the Lie groups $S$ and $K$ and choose a left-invariant Riemannian metric on this product. According to Lemma 3.1 in \cite{Pittet} and since $K$ is compact, the map $\varphi$ is bi-Lipschitz. Notice for further reference that $\widetilde{G}=\widetilde{S}\widetilde{K}$ and $Z\subseteq\widetilde{K}$ (see Theorem 5.1 and its proof in \cite{Hel}). The map 
\begin{eqnarray*}\widetilde{\varphi}:\widetilde{G}&\to&\widetilde{S}\times\widetilde{K}\\
\tilde{s}\tilde{k}&\mapsto&(\tilde{s},\tilde{k})\end{eqnarray*} 
is well-defined since $\widetilde{G}=\widetilde{S}\widetilde{K}$. Consider the following commutative diagram.
$$\xymatrix{
\widetilde{G}\ar[r]^{\widetilde{\varphi}}\ar[d]^p&\widetilde{S}\times\widetilde{K}\ar[d]^{p_1}\\
G\ar[r]^{\varphi}& S\times K
}$$
where $p_1$ is the product of the $Z$-regular cover $\widetilde{K}\to K$ with the trivial cover $\widetilde{S}\to S$. On $\widetilde{G}$ we choose the left-invariant Riemannian metric which turns $p$ into a local isometry. On $\widetilde{S}\times\widetilde{K}$, we choose the left-invariant metric (for the product structure) which turns $p_1$ into a local isometry. As $\widetilde{\varphi}$ covers $\varphi$, it is also bi-Lipschitz. Since $Z\subseteq\widetilde{K}$ is cocompact, it is undistorted, and since the inclusion $\widetilde{K}\subseteq\widetilde{S}\times\widetilde{K}$ is totally geodesic it is undistorted as well, and we conclude because $\widetilde{\varphi}$ is bi-Lipschitz.
\end{proof}
\begin{proof}[Proof of Theorem \ref{NonB=>RD}] We will show that $G$ can be expressed as a central extension of a group with property RD by means of a polynomial growth cocycle. A group $G$ as in the theorem is of the form $G=\widetilde{G}/\Gamma$, where $\Gamma$ is a discrete subgroup of $Z(\widetilde{G})$, the center of $\widetilde{G}$. Now, $Z(\widetilde{S})$, the center of $\widetilde{S}$ is discrete in $\widetilde{S}$ (see \cite{GOV}) and hence the semisimple group $\widetilde{S}/Z(\widetilde{S})$ has trivial center. The following diagram is commutative:
$$\xymatrix{
{\widetilde{G}}\ar[dr]^{p_Z}\ar[d]_{p_{\Gamma}} & \\
{G=\widetilde{G}/\Gamma}\ar[r]_p&{\widetilde{G}/Z(\widetilde{G})}
}$$
where the bottom arrow $p:G\to \widetilde{G}/Z(\widetilde{G})$ is the quotient of $G$ by $Z(\widetilde{G})/\Gamma$. Since $Z(\widetilde{G})/\Gamma$ is central in $G$, we have a central extension to which we want to apply Proposition \ref{RDextensionsPol}. To start with, 
$$\widetilde{G}/Z(\widetilde{G})=\widetilde{S}/Z(\widetilde{S})\times\widetilde{Q}/Z(\widetilde{Q})$$ 
has property RD because it is a product of two groups with property RD (see Lemma \ref{produits} combined with Theorem \ref{principal} and Proposition \ref{moyennable}). In \cite{VD}, Varopoulos proved that any closed subgroup of a connected Lie group with polynomial volume growth is at most polynomially distorted. Combined with Lemma \ref{ChP}, it implies that the center $Z(\widetilde{G})$ is at most polynomially distorted in $\widetilde{G}$. Hence $A=Z(\widetilde{G})/\Gamma$ is at most polynomially distorted in $G$. As $\widetilde{G}/Z(\widetilde{G})$ is a connected Lie group, it is generated by any neighborhood of the identity. The projection $p$ being a principal $A$-bundle map, it admits a differentiable section on a sufficiently small neighborhood of 1. Hence Proposition \ref{CocycleEtDistorsion} applies, so we conclude using Proposition \ref{RDextensionsPol}.\end{proof}
\noindent
{\it Remark.} Combining Lemma \ref{ChP} with Proposition \ref{CocycleEtDistorsion} gives a linear bound on the growth of 2-cocycles of central extensions of a semisimple Lie group with finite center. As pointed out by M. Burger and E. Ghys, it follows from results of Dupont in \cite{Dupont} that those cocycles are actually bounded. For an elementary proof, see \cite{bornes}.
%%%%%%%%%%%%%%%%%%%%%%%%%%%%%%%%%%%%%%%%%%%%%%
\section{The structure of connected Lie groups with property RD}\label{class}
%%%%%%%%%%%%%%%%%%%%%%%%%%%%%%%%%%%%%%%%%%%%%%%%%%%
In this section we finish the proof of our main theorem. To do so, we start by explaining the terms used in part (b) of Theorem \ref{classe}. Recall that a Lie algebra is of type R if all the weights of the adjoint representation are purely imaginary. A Lie group is of type R if its associated Lie algebra is of type R. According to Guivarc'h and Jenkins, a Lie algebra is of type R if and only if the associated Lie group has polynomial volume growth (see \cite{G} and also \cite{Jen}). Thus by the Fundamental Theorem of Lie (see Theorem 2.8.2  in \cite{Var}) the statements (b) and (c) in Theorem \ref{classe} are equivalent. We now turn to (a)$\Rightarrow$(b) in the proof of Theorem \ref{classe}. This part relies on Varopoulos' work in \cite{V}. Varopoulos introduces a dichotomy among finite dimensional real Lie algebras. Namely, he divides them into B-algebras and NB-algebras. We now quote the two results of \cite{V} that are crucial for our purpose.
\begin{theo}[Varopoulos, \cite{V}]\label{Bclassification} Let ${\mathbf g}$ be a unimodular algebra. Then either ${\mathbf g}$ is a B-algebra, or ${\mathbf g}$ is the direct product ${\mathbf s}\times{\mathbf q}$, where ${\mathbf q}$ is an algebra of type R and ${\mathbf s}$ is either $\{0\}$ or semisimple.\end{theo}
A Lie group is called a B-group if its Lie algebra is a B-algebra. Those groups have the following property.
\begin{theo}[Varopoulos, \cite{V}]Let $G$ be a B-group and let $d\Phi=\phi\vec{d}x$ be a compactly supported probability measure on $G$, with continuous density $\phi$. Then there exists $c>0$ such that for any $n\in\N$ we have
\begin{equation}\label{B}\phi^{*n}(1)=O\left(\|L_{\Phi^{*n}}\|_{2\ra 2}\exp(-cn^{1/3})\right).\end{equation}\end{theo}
This theorem has an easy corollary.
\begin{cor}\label{B=>NonRD} B-groups cannot have property RD.\end{cor}
\begin{proof}Let $G$ be a B-group. For $d\Phi=\phi\vec{d}x$ as in the previous theorem and $\check{\Phi}=\Phi$, we have $L_{\Phi}^*=L_{\check{\Phi}}=L_{\Phi}$ and $\|L_{\Phi^{*2n}}\|_{2\ra 2}=\|L_{\Phi^{*n}}\|_{2\ra 2}^2$. By (\ref{B}) it follows that
$$\phi^{*2n}(1)\leq A\|L_{\Phi^{*2n}}\|_{2\ra 2}\exp (-cn^{1/3})= A\|L_{\Phi^{*n}}\|^2_{2\ra 2}\exp(-cn^{1/3}),$$
for some constant $A\geq 1$. Set $f=\phi^{*n}$, so that $F=\Phi^{*n}$ and $\phi^{*2n}(1)=\|f\|^2_2$. Now assume that $G$ has property RD, then
\begin{eqnarray*}\|L_{F}\|^2_{2\ra 2}&\leq & Cn^D\|f\|^2_2=Cn^D\phi^{*2n}(1)\leq ACn^D\|L_{\Phi^{*n}}\|^2_{2\ra 2}\exp(-cn^{1/3})\\
&=&ACn^D\|L_{F}\|^2_{2\ra 2}\exp(-cn^{1/3}),\end{eqnarray*}
where $C$ and $D$ are the constants coming from the definition of property RD and the second inequality follows from the assumption that $G$ is a B-group. We conclude that $1\leq ACn^D\exp (-cn^{1/3})$, which is a contradiction for $n$ big enough. It follows that $G$ does not have property RD.\end{proof}
We can now finish the proof of our main result. It remains to show that (a) implies (b). Let $G$ be a connected Lie group. If $G$ has property RD then it has to be unimodular according to \cite{JiSch}, and if $G$ is unimodular with property RD, then it cannot be a B-group according to Corollary \ref{B=>NonRD}. The proof of Theorem \ref{classe} is now complete.\\ 

\smallskip

Recall that a group $\Gamma$ is virtually connected if the connected component of the identity is of finite index in $\Gamma$. Equivalently, the group $\Gamma$ has finitely many connected components. Recall also that any virtually connected group $\Gamma$ admits a descending sequence of compact normal subgroups $K_{\alpha}$ with $\bigcap_{\alpha}K_{\alpha}=\{e\}$ and $\Gamma/K_{\alpha}$ is a Lie group (see \cite{Glu}). We can now give a complete classification for virtually connected compactly generated groups (given by Theorem \ref{classe} combined with Lemmas \ref{extfin} and \ref{ExtCompacte}).
\begin{cor}\label{MerciBekka} Let $\Gamma$ be a virtually connected compactly generated group. Let $K$ be a normal compact subgroup such that $G=\Gamma/K$ is a Lie group. Let $G_0$ be the connected component of the identity in $G$. The following are equivalent.
\begin{itemize}
\item[(a)] $\Gamma$ has property RD.
\item[(b)] The Lie algebra ${\mathbf g}$ of $G$ decomposes as a direct product ${\mathbf g}={\mathbf s}\times{\mathbf q}$, where ${\mathbf s}$ is semisimple or $\{0\}$ and ${\mathbf q}$ is an algebra of type R.
\item[(c)] The universal cover $\widetilde{G}_0$ of $G_0$ decomposes as a direct product $\widetilde{S}\times\widetilde{Q}$, where $\widetilde{S}$ is semisimple and $\widetilde{Q}$ has polynomial volume growth.\end{itemize}\end{cor}
\bigbreak
\noindent
{\it Question.} Does Theorem \ref{classe} apply to Lie groups over $\Q_p$?\\
\bigbreak
\noindent
{\it Acknowledgements.} We thank Bachir Bekka and Alain Valette for helpful conversations and for pointing out to us Corollary \ref{MerciBekka}.
%%%%%%%%%%%%%%%%%%%%%%%%%%%%%%%%%%%%%%%%%%%%%%%%%%%%%%%%%%%%%%%%


\begin{thebibliography}{00}
%
\bibitem{Bourb} Bourbaki, N. \emph{Espaces vectoriels topologiques.} Chapitres 1 \`a 5., El\'ements de math\'ematiques.
%
\bibitem{BouIntegration} Bourbaki, N. \emph{Int\'egration.} El\'ements de math\'ematiques, livre VI.
%
\bibitem{moi}Chatterji, I. \emph{Property RD for cocompact lattices in a finite product of rank one Lie groups with some rank two Lie groups.} Geometria Dedicata {\bf 96}  (2003), 161--177.
%
\bibitem{heatker}Chatterji, Pittet and Saloff-Coste. {\em Heat decay and property RD.} In preparation.
%
\bibitem{bornes}Chatterji, Pittet and Saloff-Coste.  {\em Another proof of a result by Dupont.} In preparation.
%
\bibitem{Cow}Cowling, M. {\em Herz's "principe de majoration" and the Kunze-Stein phenomenon.}  Harmonic analysis and number theory (Montreal, PQ, 1996),  73--88, CMS Conf. Proc., {\bf 21}, Amer. Math. Soc., Providence, RI, 1997. 
%
\bibitem{CGHM} Cowling M.,  Giulini S., Hulanicki A. and Mauceri G. {\em Spectral multipliers for a distinguished Laplacian on certain groups of exponential growth.} Studia Math. {\bf 111} (1994), 103--121. 
%
\bibitem{CHH} Cowling M., Haagerup U. and Howe R. {\em Almost $L\sp 2$ matrix coefficients.} J. Reine Angew. Math. {\bf 387} (1988), 97--110.
%
\bibitem{Dupont}Dupont, J. L. \emph{Bounds for characteristic numbers of flat bundles.} Algebraic topology, Aarhus 1978, pp. 109--119, Lecture Notes in Math. {\bf 763} Springer, Berlin, 1979.
%
\bibitem{Glu}Glushkov, V.M. \emph{The structure of locally compact groups and Hilbert's fifth problem.}  Amer. Math. Soc. Transl. (2)  {\bf 15} (1960), 55--93.
%
\bibitem{G} Guivarc'h, Y. {\em Croissance polynomiale et p\'eriodes des fonctions harmoniques.} Bull. Soc. Math. France {\bf 101} (1973), 333--379. 
%
\bibitem{GOV} Gorbatsevich V.V., Onishchick, A.L. and Vinberg, E.B. {\em Foundations of Lie Theory and Lie Transformation Groups.} Springer 1997.
%
\bibitem{GNR}Gromov, M. \emph{Asymptotic invariants of infinite groups.} Proceedings of the symposium held in Sussex 1991, G. Niblo and L. Reeves Ed.
%
\bibitem{Haa}Haagerup, U. \emph{An example of nonnuclear $C*$-algebra which has the metric approximation property.} Inv. Math. {\bf 50} (1979), 279--293.
%
\bibitem{Hel}Helgason, S. \emph{ Differential geometry, Lie groups, and symmetric spaces.} Pure and Applied Mathematics {\bf 80}. Academic Press, Inc., New York-London, 1978.
%
\bibitem{HR}Hewitt, E. and Ross, K. {\em Abstract Harmonic Analysis, Vol I}. Springer, 1963
%
\bibitem{Jen} Jenkins, J. W. \emph{Growth of connected locally compact groups.}  J. Functional Analysis  {\bf 12}  (1973), 113--127.
%
\bibitem{JiSch}Ji R. and Schweitzer L. B. {\em Spectral invariance of Smooth Crossed products, and Rapid Decay for Locally compact groups.} Topology {\bf 10} (1996), 283-305.
%
\bibitem{Jol}Jolissaint, P. \emph{Rapidly decreasing functions in reduced C*-algebras of groups.} Trans. Amer. Math. Soc. {\bf 317} (1990), 167-196.
%
\bibitem{LV} Lafforgue, V. Private communication.
%
\bibitem{L1} Lafforgue, V. \emph{A proof of property RD for discrete cocompact subgroups of $SL_3({\bf R})$.} Journal of Lie Theory {\bf 10} (2000), 255--277.
%
\bibitem{Laff} Lafforgue, V. \emph{KK-th\'eorie bivariante pour les alg\`ebres de Banach et conjecture de Baum-Connes.} Invent. Math. {\bf 149}  (2002),  no. 1, 1--95.
%
\bibitem{L} Leptin, H. \emph{On locally compact groups with invariant means.} Proc. Amer. Math. Soc. {\bf 19} (1968), 489--494.
%
\bibitem{M} Mustapha, S. {\em Multiplicateurs spectraux sur certains groupes non-unimodulaires.}  Harmonic analysis and number theory (Montreal, PQ, 1996),  11--30, CMS Conf. Proc., {\bf 21}, Amer. Math. Soc., Providence, RI, 1997.
%
\bibitem{OV}Onishchik, A.L. and Vinberg, E.B. {\em Lie groups and algebraic groups.} Translated from the Russian and with a preface by D. A. Leites. Springer Series in Soviet Mathematics. Springer-Verlag, Berlin, 1990.
%
\bibitem{Pittet} Pittet, Ch. \emph{The isoperimetric profile of homogeneous Riemannian manifolds.} J. Differential Geometry {\bf 54} (2000), 255-302.
%
\bibitem{RRS}Ramagge, J., Robertson, G. and Steger, T. \emph{A Haagerup inequality for $\tilde{A_1}\times\tilde{A_1}$ and $\tilde{A_2}$ buildings.} GAFA, Vol. {\bf 8} (1988), 702--731.
%
\bibitem{Rudin} Rudin W. {\em Functional analysis.} International Series in Pure and Applied Mathematics. McGraw-Hill, Inc., New York, 1991.
%
\bibitem{WS2}Saloff-Coste L. and Woess W. {\em Transition operators, groups, norms, and spectral radii.} Pacific J. Math. {\bf 180} (1997), 333-367.
%
\bibitem{WS3} Saloff-Coste L. and Woess W. {\em Transition operators on co-compact $G$-spaces.} Preprint.
%
\bibitem{Val}Valette, A. \emph{Introduction to the Baum-Connes Conjecture.} Lectures in Mathematics ETH Z\"urich, Birkh\"auser Verlag, Basel, 2002.
%
\bibitem{Var} Varadarajan, V.S. {\em Lie groups, Lie algebras, and their representations.} Graduate Texts in Mathematics, Springer-Verlag, 1984.
%
\bibitem{VSC} Varopoulos, Saloff-Coste and Coulhon. {\em Analysis and geometry on groups.} Cambridge University Press, 1993.
%
\bibitem{VD} Varopoulos, N. \emph{Distance distortion on Lie groups.} Random walks and discrete potential theory (Cortona, 1997), 320--357,Sympos. Math., Cambridge Univ. Press, Cambridge, 1999. 
%
\bibitem{V} Varopoulos, N. {\em Analysis on Lie groups.} Revista Matematica Iberoamericana Vol 12, No3, 1996. 
%

\bibitem{Woess} Woess, W. {\em  Random walks on infinite graphs and groups.} Cambridge Tracts in Mathematics, {\bf 138}. Cambridge University Press, Cambridge, 2000. 
\end{thebibliography}
\end{document}